# Gautchi's ratio and the Volume of the Unit Ball in $R^n$

Dimitris Karayannakis,
Department of Science / Section of Mathematics,
T.E.I. of Crete, Heraklion 71004,Greece
**dkar@stafff.teicrete.gr**

Abstract

Let $\Omega_n$ be the volume of the unit ball in $R^n$. We formulate as an infinite product the gamma function ratio $\Gamma(x+1/2)/\Gamma(x), x > 0$, which allows us to reproduce and /or produce a variety of formulas and inequalities, some of the latter seemingly new, concerning $\Omega_n$, $\Omega_{n-1}/\Omega_n$ and $\Omega_n^2/\Omega_{n-1}\Omega_{n+1}$.

**1. Introduction**

Using the gamma function $\Gamma(x) = \int_0^\infty u^{x-1}e^{-u}du$, $x > 0$, we can represent the volume of the $n$-dimensional Euclidean ball, $\Omega_n = \pi^{\frac{n}{2}}/\Gamma(1+\frac{n}{2})$ for any positive integer $n$ and with the convention $\Omega_0 = 1$ this volume formula holds for $n \geq 0$. In addition the so called psi function $\psi(x)$ is by definition the logarithmic derivative of the gamma function, $\Gamma'(x)/\Gamma(x)$ $x \notin \mathbb{Z}_0^-$.
A basic functional equation concerning the psi (or digamma) function, that can be easily derived from gamma's functional equation $\Gamma(x+1) = x\Gamma(x)$, reads $\psi(x+1) = \psi(x) + 1/x$ and a well known series repesentation of psi is:

$$\psi(x+1) = -\gamma + \sum_{k=1}^\infty (\frac{1}{k} - \frac{1}{k+x}), \text{ for } x > -1 \qquad (1.1)$$

where $\gamma = 0.577215..$ (the Euler-Mascheroni constant) and $\psi(1/2) = -\gamma - \log 4$.
We will also make use of the following classical formula:

$$\psi(u) - \psi(v) = \sum_{k=0}^\infty (\frac{1}{k+v} - \frac{1}{k+u}), \text{ for } u \notin \mathbb{Z}_0^-, v \notin \mathbb{Z}_0^- \qquad (1.2)$$

(A comprehensive presentation of these two special functions can be found in Chapter 1 of [3]. For $\Omega_n$ see [2] and/or [3] and the referencew within.)

In [6] (Lemma II.1) using what is probably the oldest definition of the gamma function from Euler's correspondence in 1729), namely

$$x\Gamma(x) = \prod_{k=1}^{\infty} \frac{k^{1-x}(k+1)^x}{x+k}, \text{ we obtained that for } x > 0, a < 1 \text{ and } x + a \notin \mathbb{Z}_0^-$$

$$\frac{\Gamma(x+a)}{\Gamma(x)} = \frac{f(x,a)}{\Gamma(1-a)} \quad \text{with} \quad f(x,a) = \prod_{k=1}^{\infty} \frac{k(x+k-1)}{(k-a)(k+x+a-1)} \tag{1.3}$$

The ratio $\Gamma(x+a)/\Gamma(x)$ is usually called Gautchi's ratio in the literature.
In [7] we have called the infinite product $f(x,a)$ in (1.3) the joint factor of $\Gamma(x+a)$ and for convenience we have denoted the product of the first $m$ terms of this infinite product as $f_m(x,a)$ which we called the $m$-truncate of the joint factor. It is also evident that this representation of Gautchi's ratio for the special case $x = 1-a$ should be considered equal to 1, a convention compatible with the infinite product definition and the fact $\Gamma(1) = 1$. Various numerical applications of this and other revisited formulas for a rational argument $x$ have been presented in [7] where a heuristic construction of this factorization has been given '"from scratch" (by formal use of power series).

In this work we present various other applications that involve specifically Gautchi's ratio for $a = 1/2$. Using the corresponding joint factor we will reformulate the sequences $(\Omega_{n-1}/\Omega_n)$ and $(\Omega_n^2/\Omega_{n-1}\Omega_{n+1})$ and then by use of their $m$-truncates, not only we will be able to recapture some relatively recent results, but also to produce, seemingly new inequalities concerning these sequences.

**2. Lower Bounds for $\Omega_{n-1}/\Omega_n$**

Let $\Omega_n$ be the volume of the unit ball in $\mathbf{R}^n$, $n \in \mathbf{N}$. Since $\Omega_n = \pi^{\frac{n}{2}}/\Gamma(1+\frac{n}{2})$, setting in (1.3) $a = \frac{1}{2}$ and $x = \frac{n+1}{2}$, we can immediately conclude that

$$v_n = \Omega_{n-1}/\Omega_n = f(\frac{n+1}{2}, \frac{1}{2})/\pi = \frac{1}{\pi}\prod_{k=1}^{\infty}(\frac{2k}{2k-1})(\frac{2k+n-1}{2k+n}). \tag{2.1}$$

It is easy to check that all the terms inside the infinite product in (2.1) exceed 1 and then we directly obtain for any $m \in \mathbf{N}$ the following inequality:

$$f_m(\frac{n+1}{2}, \frac{1}{2})/\pi < v_n \tag{2.2}$$

For simplicity we will denote the $m$ truncate in (2.2) by $f_m(n)$

The lower bounds in (2.2) seem to be the simplest rational (in $n$) bounds, holding for all $n$ and naturally the larger the value of $m$ the sharper they become. In particular for $m = 1$ we obtain the simplest and crudest among them all, namely $2(n+1)/\pi(n+2)$. It is evidently less sharp, for all $n \neq 1$, compared to $b_n = \sqrt{n/2\pi}$ presented in 1987 by Borgwardt who also derived $b_{n+1}$ as an upper bound (see the discussion in [1]).

If we still wish to exploit (2.1), without settling for an m truncation only as we did in (2.2), and also refine the derived lower bound of $\Omega_{n-1}/\Omega_n$, a switch from the infinite product formulation to its equivalent logarithmic series formula together with (1.1) will suffice. These new lower bounds that we will obtain below turn out to be less elegant compared to Borgwardt's, or compared to our m bounds in (2.2), but even the crudest of them will be sharper than $b_n$:

First observe that using the inequality $\log t \geq 1 - 1/t, t > 0$, we have for $n \neq 1$

$$\log f(\frac{n+1}{2}, \frac{1}{2}) \geq \log f_m(n) + [\frac{n}{2}(\Sigma_1 - \sigma_m)] \tag{2.3}$$

where $\Sigma_1 = \frac{1}{n-1}\sum_{k=1}^{\infty}(\frac{1}{k} - \frac{1}{k+(n-1)/2})$ and $\sigma_m = \sum_{k=1}^{m}\frac{1}{k(2k+n-1)}$.

Using (1.1) we can immediately conclude that for $n \neq 1$ and all m

$$d_n(m) = \frac{1}{\pi}f_m(n)\exp\frac{n}{2}[\frac{\psi((n+1)/2) + \gamma}{n-1} - \sigma_m] < v_n \tag{2.4}$$

For $n = 1$, a simple application of l'Hospital's rule along with the fact that $\sum_{k=1}^{\infty}1/k^2 = \psi'(1) = \pi^2/6$ (see [5]) allows us the use of (2.4) for all n.

It is worth mentioning at this point that $\psi'(x) = \sum_{k=0}^{\infty}1/(k+x)^2, x \notin \mathbb{Z}_0^-$, is called the trigamma function and we will require some additional results about this function later on in Application 4, where we seek a lower bound for the volume ratio $\Omega_n^2/\Omega_{n-1}\Omega_{n+1}$.

Despite the awsome looks of (2.4) we notice that calculating $\psi((n+1)/2)$ at least for small values of n becomes rather easy by use of psi's functional equation and although an approximation of the value of $\gamma$ will be involved in proving later on (2.6), its presence in (2.4) is practically "virtual" since it vanishes in the actual calculations due to (1.1).

The crudest and less cumbersome among the lower bounds given by (2.4) is evidently $d_n(1)$ which using an ordinary scientific pocket calculator leads to lower bounds more crude than Alzer's $A_n = \sqrt{(n+1/2)/2\pi}$ ([1]). On the other hand we will show that $d_n(1) = d_n > b_n$ or, equivalently, $\log d_n > \log b_n$.
Firstly, we can easily check from (2.4) that

$$d_n = \frac{2}{\pi}(\frac{n+1}{n+2})\exp\frac{n}{2}[\frac{\psi((n+1)/2) + \gamma}{n-1} - \frac{1}{n+1}]$$

Now due to technicalities in the proof that follows we check separately the cases $n = 1$ and $n = 2$ and we find that:

$$d_1 = \exp[(\frac{\pi^2}{6} - 1)/4] > 0.4986 > b_1 = 0.4886...,$$

and $d_2 = \frac{4}{8\pi} e^{5/3} > 0.6319 > b_2 = 0.5641...$

In order to show that $\log d_n > \log b_n$, for $n \geq 3$, will make use of the folowing inequality (Th.1 in [4]) :

$$\log(t+1/2) - 1/t < \psi(t), \text{ for } t > 0 \tag{2.5}$$

Combining $\gamma > 1/2$, with (2.5) for $t = (n+1)/2$, it will suffice to show that

$$\frac{n}{n-1}[\log(\frac{n+2}{2}) - \frac{1}{2(n-1)}] \geq \log \frac{\pi}{8} n(\frac{n+2}{n+1})^2, \text{ for } n \geq 2, \tag{2.6}$$

Since $\log 2 > 1/2$ we can easily see that in order (2.6) to be true, it suffices to show

$$\frac{1}{2(n-1)} \leq \log[\frac{4}{\pi} \frac{(n+1)^2}{n(n+2)}], \text{ for } n \geq 3, \tag{2.7}$$

and this is clearly true since the l.h.s term of (2.7) is $\leq 1/6$ while the r.h.s term is $> \log(4/\pi)$ ($= 0.24...$) for all $n$.

Remark 1

We can automatically conclude from (2.1) that

(a) There exists an $r = r(n)$ such that

$$\prod_{k=1}^{r} (\frac{2k}{2k-1})(\frac{2k+n-1}{2k+n}) > \pi A_n = \sqrt{\pi(2n+1)}/2,$$

(b) For all $m \in N$, $\prod_{k=1}^{m} (\frac{2k}{2k-1})(\frac{2k+n-1}{2k+n}) < \pi b_{n+1} = \sqrt{\pi(2n+2)}/2$

### 3. Upper Bounds for $\Omega_{n-1}/\Omega_n$

In a similar way, this time using $\log t \leq t - 1, t > 0$, we obtain the upper bounds

$$h_n(m) = \frac{1}{\pi} f_m(n) \exp[n(\Sigma_2 - s_m)] \tag{3.1}$$

with $\Sigma_2 = \frac{1}{2(n+1)}[\sum_{k=0}^{\infty}(\frac{1}{k-1/2} - \frac{1}{k+n/2})], s_m = \sum_{k=0}^{m} \frac{1}{(2k-1)(2k+n)}.$ (3.2)

Using (1.2) we can immediately conclude that for all $n$ and $m$

$$h_n(m) = \frac{1}{\pi} f_m(n) \exp n[\frac{\psi(n/2) - \alpha}{2(n+1)} - s_m] > v_n \tag{3.3}$$

where $\alpha = \psi(-1/2) = \psi(3/2) = 2 - \gamma - \log 4 \, (\cong 0.036)$ (3.4)

Remark 2

(a) Using graphing via Mathematica we can check that $h_n(1)$ is less sharp, for all $n \in N$, than Alzer's upper bound $\left(\dfrac{n+\beta}{2\pi}\right)^{1/2}$ with $\beta = \pi/2 - 1$, given in [1], which for $n = 1$ gives the exact value of $v_1$. We can also check that even

$$h_n(2) = \frac{8}{3\pi} \frac{(n+1)(n+3)}{(n+2)(n+4)} \exp[n\frac{\psi(n/2) - \alpha}{2(n+1)} - \frac{n^2 - 4n - 24}{3(n+1)(n+2)(n+4)}]$$

is sharper than Alzer's only for $2 \leq n \leq 7$ e. t. c. Thus we oberve that, for $m \geq 2$, there exists an increasing sequence of natural numbers $r_m > 2$ such that $h_n(m) < \left(\dfrac{n+\beta}{2\pi}\right)^{1/2}$ if and only if $2 \leq n \leq r_m$.

(b) Similarly, we can check that $h_n(1) = \dfrac{2}{\pi}(\dfrac{n+1}{n+2}) \exp[n\dfrac{\psi(n/2) - \alpha}{2(n+1)} + \dfrac{2}{n+2}]$

is sharper than Bogwardt's $\left(\dfrac{n+1}{2\pi}\right)^{1/2}$ only for $1 \leq n \leq 7$, resulting thus in an analogous claim to (a) that will include this time also the case $m = 1$.

**4. Lower Bounds for $\Omega_n^2 / \Omega_{n-1}\Omega_{n+1}$**

Using (2.1) we can immediate conclude that for all $n \in N$

$$w_n = \Omega_n^2 / \Omega_{n-1}\Omega_{n+1} = v_{n+1}/v_n = \prod_{k=1}^{\infty} \frac{(2k+n)^2}{(2k+n)^2 - 1} \qquad (4.1)$$

Once again, all the terms inside the infinite product in (4) exceed 1 and then we directly obtain, for any $m \in N$, the following inequality:

$$\prod_{k=1}^{m} \frac{(2k+n)^2}{(2k+n)^2 - 1} < w_n \qquad (4.2)$$

and the larger the value of $m$ the sharper these rational lower bounds become. If we wish to exploit (4.1) without resorting to (4.2) for a big $m$, that would had led to an extremely cumbersome lower bound, we can again appeal to the logarithm of (4.1). By use once more of the inequality $\log t > 1 - 1/t, t \neq 1,, t > 0$, we obtain the lower bound $\exp(\dfrac{1}{4}\psi'(\dfrac{n}{2}) - \dfrac{1}{n^2})$.

Then by use of the inequality (see (26) in [8])

$$\psi'(x) > \frac{1}{x+1} + \frac{1}{x^2} + \frac{1}{2(x+1)^2}, \text{ for } x > 0, \qquad (4.3)$$

we obtain for all $n \in \mathbb{N}$ the inequality $\exp[(n+3)/2(n+2)^2] < V_n$     (4.4)
As we can easily check, for large $n$, the lower bounds in (4.4) are still not so satisfactory. But if we isolate the leading term in the infinite product of (4.1) and work in a similar way as we did before we obtain the lower bound

$$\frac{(n+2)^2}{(n+1)(n+3)} \exp\left(\frac{1}{4}\psi'\left(\frac{n}{2}\right) - \frac{1}{n^2} - \frac{1}{(n+2)^2}\right)$$ and using (4.3) we obtain

$$p_n = \frac{(n+2)^2}{(n+1)(n+3)} \exp[(n+1)/2(n+2)^2] < w_n. \quad (4.5)$$

Remark 3

(a) An interesting point concerning (4.2) is that even $(n+2)^2/(n+1)(n+3)$, the crudest (but most elegant) among all these lower bounds exceeds 1, thus re-establishing the known inequality $1 < w_n$ without involving (directly) the pertinent monotonicity theory and the logarithmic concavity of the sequence $\Omega_n, (n=0,1,...)$ (see the discussion in [1] or [2]).

(b) The bound in (4.5) consists a numerical improvement not only of the bound in (4.4) but also of the known lower bound $(1+\frac{1}{n+1})^{1/2}$ (see the discussion in [1] or [2]). To see this we must show that $2\log p_n > (n+2)/(n+1)$ or equivalently that $\frac{n+1}{(n+2)^2} > \log\frac{(n+1)(n+3)^2}{(n+2)^3}$.

Thus it will suffice to show that $\frac{n+1}{(n+2)^2} \geq \frac{(n+1)(n+3)^2}{(n+2)^3} - 1$,

or $\frac{n+1}{(n+2)^2} \geq \frac{n^2+3n+1}{(n+2)^3}$, which is evidently true.

**5.** An Upper Bound for $\Omega_n^2/\Omega_{n-1}\Omega_{n+1}$

We derive now an upper bound for $w_n$, using again the pertinent classical logarithmic inequality along with (1.2) and the functional equation of psi. From (4.1) and for $n \neq 1$ we obtain that

$$\log w_n < \frac{1}{4}\left[\psi\left(\frac{n+1}{2}\right) - \psi\left(\frac{n-1}{2}\right)\right] - \frac{1}{n^2-1} = \frac{1}{4}\left(\frac{2}{n-1}\right) - \frac{1}{n^2-1} = \frac{0.5}{n+1}.$$

Thus $w_n < \exp[0.5/(n+1)]$     (5.1)

If $n = 1$ we have $\log w_1 < \frac{1}{4}\sum_{k=1}^{\infty}\frac{1}{k(k+1)} = \frac{1}{4}$ and so (5.1) holds for all $n$.

Remark 4

(a) The bound in (5.1) not only improves Alzer's $(1+1/n)^{1/2}$ but also improves M. Merkle's upper bound $(n+2)^{3/2}/(n+1)(n+3)^{1/2}$ (Th.5 in [8]) which in its turn constitutes an improvement of Alzer's.

Note that this should have been the correct form of the r.h.s of (42) in [8] instead of the printed one). M.Merkle pointed out ( personal communication) that a better improvement (for $n \geq 2$) can be obtained by use of (23) in [8] instead of (22) and then it would had given the upper bound $(1+2/n)^{1/4}$.

(b) We show now that our upper bound in (5.1) also refines this upper bound. Equivalently it suffices to show that $2/(n+1) < \log(1+2/n)$: for $x > 0$ we set $f(x) = \log(1+2/x) - 2/(x+1)$. Then

$f'(x) = 2[1/x(x+2) - 1/(x+1)^2] > 0$, and so $f(x)$ is strictly monotonically increasing. Since $f(1) = \log 3 - 1 > 0$ we conclude that $f(n) > 0$ for all $n$.

(c) Since (5.1) immediately implies $w_n < 1 + 1/n$ we have re-established by classical analysis means, the strictly increasing nature of the sequence $n\Omega_n / \Omega_{n-1}$, as it was first proved in 1997 by Klein and Rota.

(d) The bound in (5.1) can be improved even further if we combine the $m$-truncate of (4.1) with the tail of the series that led to (5.1). For example with $m = 1$ we obtain the sharper upper bound $\dfrac{(n+2)^2}{(n+1)(n+3)} \exp[0.5/(n+3)]$.